\documentclass{amsart} \usepackage[all]{xy} \usepackage[latin1]{inputenc} \usepackage{amsmath,amsthm,pb-diagram,pb-xy,amscd, amssymb, amsfonts}
\usepackage[dvips]{graphicx} \usepackage{hyperref}
\numberwithin{equation}{section}

\newcommand{\bb}{\mathcal{B}}

\newcommand{\End}{\mbox{\rm End\,}}
\theoremstyle{plain}

\newcommand{\Hom}{\text{Hom}}
\newcommand{\Aut}{\text{Aut}}

\newcommand{\Obj}{\text{Obj}}

\newcommand{\C}{\mathcal{C}}

\newcommand{\D}{\mathcal{D}}
\newcommand{\M}{\mathcal{M}}

\newcommand{\id}{\mathrm{id}}

\newtheorem{teor}{Theorem}[section] 
\newtheorem{ejem-prop}[teor]{Example} \newtheorem{prop}[teor]{Proposition}
\newtheorem{lem}[teor]{Lemma}

\theoremstyle{definition} \newtheorem{defin}[teor]{Definition} 
\newtheorem{ejem}[teor]{Example}

\theoremstyle{remark} \newtheorem{obs}[teor]{Remark}

\begin{document}

\title[Crossed product tensor categories]{Crossed product tensor categories} \author{C\' esar
Galindo } \address{Departamento de Matem?ticas\newline \indent
Universidad de los Andes \newline \indent Bogot?, Colombia}
\email{cn.galindo1116@uniandes.edu.co, cesarneyit@gmail.com}
 \date{\today}

\begin{abstract} A graded tensor category over a group $G$ will be called a  crossed product tensor category if every homogeneous component has at least one multiplicatively invertible object. Our main result is a description of crossed product tensor categories, graded monoidal functors, monoidal natural transformations,  and braidings in terms of coherent outer $G$-actions over tensor categories. \end{abstract}

\maketitle

\section{Introduction}

A $G$-graded ring $A= \bigoplus_{\sigma\in G}A_\sigma$ is called a $G$-crossed product if each $A_\sigma$ has an invertible element.  Some important classes of rings like skew group-rings and  twisted group-rings are special cases of crossed product rings. One of the basic examples is the group algebra of a group $F$, it is graded by a quotient group of $F$, see  \cite[Subsection 11C]{CR}. In this case, the representation theory of $F$ can be analyzed using Clifford theory, see \cite[Subsection 11A]{CR}.
\smallbreak
In analogy with graded rings, a $G$-graded tensor category (see Subsection \ref{defincion graded an crossed}) $\C=\bigoplus_{\sigma\in G}\C_\sigma$ will be called $G$-crossed product tensor category, if there is an invertible object in each homogeneous component of $\C$.
\smallbreak
Graded tensor categories appear naturally in classification problems of fusion categories and finite tensor categories \cite{ENO}, \cite{finite-categories}. One of the most interesting examples of $G$-crossed product tensor categories is the semi-direct product tensor category associated to an action of a group over a tensor category, see \cite{tambara}. Semi-direct tensor product categories have been used in order to solve an important open problem in semisimple Hopf algebras theory \cite{nik}.
\smallbreak

Crossed product rings are commonly described using crossed systems \cite{MoGR}. Crossed systems can be interpreted in terms of monoidal functors in the following way: if $A$ is a ring, let we denote by $\overline{Out(A)}$ the monoidal category (in fact, it is a categorical-group) of outer automorphisms, where the objects are automorphisms of $A$, and the arrows between automorphisms $\sigma$ and $\tau$ are  invertible elements $a\in A$, such $a\sigma(x)= \tau(x)a$ for all $x\in A$. Given a group $G$, the data that define $G$-crossed systems are the same as the data that define  monoidal functors from $\underline{G}$ to $\overline{Out(A)}$, where $\underline{G}$ is the discrete monoidal category associated with $G$.
\smallbreak
We develop  crossed product system theory for crossed product tensor categories using higher category theory. To do so, we use the monoidal bicategory $\text{Bieq}(\C)$ (in fact, it is a weak 3-group), whose 0-cells are tensor autoequivalence of $\C$, 1-cells are pseudonatural isomorpshism, and 2-cells are modification (see Subsection \ref{Bieq(C)}). So, a $G$-crossed system or coherent outer $G$-action over $\C$ must be a trihomomorphism from  $\underline{\underline{G}}$ (the discrete 3-category associated to $G$, see \emph{Remark} \ref{3-group discreto}) to $\text{Bieq}(\C)$.
\smallbreak
The main goal of this paper is to describe the 2-category of $G$-crossed pro\-duct tensor categories in terms of coherent outer $G$-actions over a tensor category (Theorem \ref{main theorem}), and describe the braidings of $G$-graded tensor categories (Theorem \ref{teor braid}).
\smallbreak

The main motivation for this work was the paper \cite{G}. In \textit{loc. cit.} was proposed a Clifford theory categorification for crossed product tensor categories, in order to describe simple module categories in terms of subgroups and induced module categories. We stress that in \cite{G} crossed product tensor categories were called strongly graded tensor categories.
\smallbreak

While this paper was at final stages of preparation,
Etingof, Nikshych, and Ostrik posted the preprint \cite{ENO3} containing results similar to some of ours. In \cite{ENO3} they study fusion categories graded by a group, using invertible bimodule categories over fusion categories. They reduce the classification problem of fusion categories graded by a group $G$ to classification (up to homotopy) of maps from $BG$ to classifying spaces of certain higher groupoids. In  Section \ref{seccion homotopy} we explain briefly the connection of our results with some results in \cite{ENO3}.
\smallbreak

The organization of the paper is as follows: in Section \ref{preliminares} we recall the main definitions of bicategory theory, as well as the definitions of categorical-groups, graded tensor categories, and the monoidal structure over $\text{Bieq}(\C)$. In Section \ref{section action y correspondecia} we define incohe\-rent and coherent outer $G$-actions over a tensor category $\C$, and we show an explicit bijective correspondence between equivalence classes of $G$-crossed product tensor categories and coherent outer $G$-actions. In Section \ref{section classif} we show a biequivalence between the 2-category of crossed product tensor categories and the 2-category of coherent outer $G$-actions. In Section \ref{seccion 6} we provide necessary and sufficient conditions for the existence of a braiding over a crossed product tensor category. Finally, in Section \ref{seccion homotopy} we explain the connection of our results with some results in \cite{ENO3}

\section{preliminaries}\label{preliminares}

\subsection{General conventions}
Throughout this article we work over an arbitrary field $k$. By  a tensor category $(\C,\otimes,\alpha, I )$ we understand a $k$-linear abelian category $\C$,
endowed with a $k$-bilinear exact bifunctor $\otimes:\C\times \C\to \C$, an object $I\in \C$, and
an associativity constraint  $\alpha_{V,W,Z}: (V\otimes W)\otimes Z\to V\otimes (W\otimes Z)$, such that
Mac Lane's pentagon axiom holds \cite{Kassel},   $V\otimes I= I\otimes V = V$, $\alpha_{V, I,W}= \text{id}_{V\otimes
W}$ for all $V, W\in \C$ and $\dim_k \End_\C(I)=1$.
\smallbreak
We shall consider only monoidal categories in which the unit constraints are identities. So, without loss of generality, we shall suppose that for every monoidal functor  $(F,\psi):\C\to \D$, we have $F(I_\C)=I_\D$ and $\psi_{V,I}=\psi_{I,V}=\id_V$, since each monoidal functor is  monoidally equivalent to one with these properties.

\subsection{Bicategories}\label{seccion bicategorias}

In this section we review some definitions on bicategory theory that we
shall need later. We refer the reader to \cite{Be} for a detailed exposition on the
subject.

\begin{defin}
A bicategory  $\mathcal{B}$ consists of the following data
\begin{itemize}
  \item a set $\Obj(\mathcal{B})$ (with elements $A,B,\ldots $ called 0-cells ),
  \item for each pair $A, B \in \Obj(\mathcal{B})$, a category $\mathcal{B}(A,B)$ (with objects  $V, W, \ldots$ called 1-cells and morphisms $f, g, \ldots$ called 2-cells),
  \item for each  $A, B, C\in \Obj(\mathcal{B})$, a bifunctor $$\overline{\circ}^{ABC}:\mathcal{B}(A,B)\times \mathcal{B}(B,C)\to \mathcal{B}(A,C),$$
  \item for each 0-cell $A\in \Obj(\mathcal{B})$, a 1-cell $I_A\in \mathcal{B}(A,A)$,
  \item for each $A, B,C, D \in \Obj(\mathcal{B})$,  natural isomorphisms (constraint of associativity)
  \begin{align*}
    \alpha^{A,B,C,D}:  - \overline{\circ}^{ABD} ( - \overline{\circ}^{BCD} - ) &\to ( - \overline{\circ}^{ABC} - ) \overline{\circ}^{ACD} -:\\
    \mathcal{B}(A,B)\times \mathcal{B}(B,C)\times \mathcal{B}(C,D)&\to \mathcal{B}(A,D).
    \end{align*}
\end{itemize}
Subject to the following axioms
\begin{itemize}
  \item coherence of the associativity:
if $(S,T,U,V)$ is an object in $\mathcal{B}(A,B)\times
\mathcal{B}(B,C)\times \mathcal{B}(C,D)\times \mathcal{B}(D,E)$,
the next diagram commutes
$$
\begin{diagram}
\node{}\node{S\overline{\circ}(T\overline{\circ}(U\overline{\circ} V))} \arrow{sw,l}{\alpha_{S,T,U\overline{\circ} V}}\arrow{se,l}{\id\overline{\circ} \alpha_{T,U,V}}\node{}\\
\node{(S\overline{\circ} T)\overline{\circ} (U\overline{\circ} V)}\arrow{s,r}{\alpha_{S\overline{\circ} T,U,V}} \node{}\node{S\overline{\circ}((T\overline{\circ} U)\overline{\circ} V)} \arrow{s,r}{\alpha_{S,T\overline{\circ} U,V}}\\
\node{((S\overline{\circ} T)\overline{\circ} U)\overline{\circ} V} \node{} \node{(S\overline{\circ} (T\overline{\circ} U))\overline{\circ} V} \arrow[2]{w,t}{\alpha_{S,T,U}\overline{\circ} \id}
\end{diagram}
$$
\item coherence of the unity
\[
\alpha_{S,I_B,T}=\id_{S\overline{\circ} T}.
\]
\end{itemize}
\end{defin}
If $\alpha$ is the identity, we have $(S\overline{\circ} T)\overline{\circ} U = S\overline{\circ} (T\overline{\circ} U)$ and similarly for morphisms, in this case we shall say that $\mathcal{B}$ is a 2-category.
\smallbreak
A monoidal category $(\C, \otimes, \alpha, I)$ is the same as a bicategory with only one 0-cell, and in this case $\overline{\circ}=\otimes$.
\begin{defin}
Let $\mathcal{B}=(\overline{\circ}, I,\alpha)$ and
$\mathcal{B}'=(\overline{\circ}, I',\alpha')$ be bicategories.
A pseudofunctor  $\Phi=(F,\phi)$ from $\mathcal{B}$ to $\mathcal{B}'$
consists of the following data
\begin{itemize}
  \item a function $F: \Obj(\mathcal{B})\to \Obj(\mathcal{B}')$, $A\mapsto F(A)$,
  \item for each pair $A,B\in \Obj(\mathcal{B})$, functors  $$F_{AB}:\mathcal{B}(A,B)\to \mathcal{B}'(F(A),F(B)),\   \ S\mapsto F(S), \  f\mapsto F(f),$$
  \item for each triple $A, B, C  \in \Obj(\mathcal{B})$, a natural isomorphism  $$\psi^{ABC}:  F_{AC}(- \overline{\circ}^{ABC} -)\to F_{AB}(-)\overline{\circ}^{F(A)F(B)F(C)} F_{BC}(-),$$
\end{itemize}
Subject to the following axioms
\begin{itemize}
  \item[(i)] $F_{AA}(I_A)= I'_{F(A)}$,
  \item[(ii)] if $(S,T,U)$ is an object in $\mathcal{B}(A,B)\times \mathcal{B}(B,C)\times \mathcal{B}(C,D)$, the following diagram commutes (where the indexes have been omitted)
\[
\begin{diagram}
\node{F(S\overline{\circ} (T\overline{\circ} U))}\arrow{s,l}{\psi}\arrow[3]{e,t}{F(\alpha)}\node[2]{}\node{F((S\overline{\circ} T)\overline{\circ} U )}\arrow{s,r}{ \psi}\\
\node{F(S)\overline{\circ} F(T\overline{\circ} U)}\arrow{s,l}{\id\overline{\circ}\psi}\node[2]{}\node{F(S\overline{\circ} T)\overline{\circ} F(U)}\arrow{s,r}{\psi\overline{\circ} \id}\\
\node{F(S)\overline{\circ} (F(T)\overline{\circ} F(U))}\arrow[3]{e,t}{\alpha'}\node[2]{}\node{(F(S)\overline{\circ} F(T))\overline{\circ} F(U)}\\
\end{diagram}
\]
  \item[(iii)] if $S$ is an object in $\mathcal{B}(A,B)$, then $\psi_{S,I_B} =\id_{F(S)}$ and $\psi_{I_A,S}= F(S)$,
  for each pair of 0-cells $A, B \in \Obj(\mathcal{B})$.
\end{itemize} \end{defin}

\begin{obs}
\begin{enumerate}
  \item The notion of pseudofunctor can be in some manner dualized by reversing the direction of the 2-cells $F_{A,B}$, this notion will be called op-pseudofunctor.
  \item A pseudofunctor between monoidal categories is just a monoidal functor and an op-pseudofunctor is an op-monoidal functor.
\end{enumerate}
\end{obs}

\begin{defin}
Let   $F,G:\mathcal{B}_0\to \mathcal{B}_1$ be  pseudofunctors between
bicategories $\mathcal{B}_0$ and $\mathcal{B}_1$. A pseudonatural transformation  $\sigma: F\to G$,
consists of the following data
\begin{itemize}
  \item  for each $A\in \Obj(\mathcal{B}_0)$,  1-cells $\sigma_A\in \mathcal{B}_1(F(A),G(A))$,
  \item for each pair $A, B\in \Obj(\mathcal{B}_0)$, and each 1-cell $V\in \mathcal{B}_0(A,B)$ a natural isomorphism
$$\sigma_V:F^{AB}(V)\overline{\circ}^{F(A)F(B)G(B)} \sigma_B \to \sigma_A \overline{\circ}^{F(A)G(A)G(B)} G^{AB}(V),$$
\end{itemize}
such that $\sigma_{I_A}=\id_{I_A}$ for all $A\in \Obj(\mathcal{B}_0)$ and for all $S\in \mathcal{B}_0(A,B)$, $T\in
\mathcal{B}_0(B,C)$, the following diagram
$$
\begin{diagram}
\node{ F(S T)\sigma}\arrow{s,l}{ \psi^F\overline{\circ} \id_\sigma}\arrow[2]{e,t}{\sigma_{ST}}\node[2]{\sigma G(S T) }\arrow{s,r}{\id_\sigma\overline{\circ} \psi^{G}}\\
\node{ F(S)F(T)\sigma}\arrow{e,t}{\id_{F(S)}\overline{\circ}\sigma_T} \node{F(S) \sigma G(T)}\arrow{e,t}{ \sigma_S\overline{\circ} \id_{G(T)}}\node{\sigma G(S) G(T)}
\end{diagram}
$$ commutes (where associativity  constraint, indexes, and the symbols $\overline{\circ}$ between objects have been omitted as a space-saving measure).
\smallbreak
\end{defin}
\begin{obs}
Again, the notion of pseudonatural transformation can be dualized by reversing the order of the natural isomorphisms $\sigma_V$, this notion will be called op-pseudonatural transformation.
\end{obs}

Pseudonatural transformations may be composed in the obvious way. If $\sigma:F\to G$, and $\tau: G\to H$ are pseudonatural transformations, then we define a new pseudonatural transformation  $\sigma\overline{\circ}\tau:F\to H$ by $(\sigma\overline{\circ}\tau)_A=\sigma_A\overline{\circ}\tau_A$, and $(\sigma\overline{\circ}\tau)_V$ is defined by the commutativity of the diagram:

\[
\begin{diagram}
\node{F(V)\overline{\circ}(\sigma_A\overline{\circ}\tau_A)} \arrow{s,l}{\alpha_{F(V),\sigma_A,\tau_A}} \arrow[3]{e,t}{(\sigma\overline{\circ}\tau)_V}\node[2]{} \node{(\sigma_B\overline{\circ}\tau_B)\overline{\circ} H(V)}\\
\node{(F(V)\overline{\circ}\sigma_A)\overline{\circ}\tau_A} \arrow{s,l}{\sigma_V\overline{\circ}\id_{\tau_A}}\node[2]{} \node{\sigma_B\overline{\circ}(\tau_B\overline{\circ} H(V))}\arrow{n,r}{\alpha_{\sigma_B,\tau_B,H(V)}}\\
\node{(\sigma_B\overline{\circ} G(V))\overline{\circ}\tau_A}\arrow[3]{e,t}{\alpha^{-1}_{\sigma_B,G(V),\tau_A}} \node[2]{}\node{\sigma_B\overline{\circ}(G(V)\overline{\circ}\tau_A)} \arrow{n,r}{\id_{\sigma_B}\overline{\circ}\tau_V}\\
\end{diagram}
\]
where the index have been omitted.

A modification between two pseudonatural transformations $\Gamma:\sigma\to \widetilde{\sigma} $, consists of 2-cells $\Gamma_A: \sigma_A \to \widetilde{\sigma}_A$   in $\mathcal{B}_1(F(A),G(A))$, such that  for all 1-cell $V\in \mathcal{B}_0(A,B)$ the diagram
$$
\begin{diagram}
\node{F^{AB}(V)\overline{\circ} \sigma_B}\arrow{s,l}{ \text{id}\overline{\circ} \Gamma_B}\arrow{e,t}{\sigma_V}\node{\sigma_A \overline{\circ} G^{AB}(V)} \arrow{s,r}{ \Gamma_A\overline{\circ}\id}\\
\node{ F^{AB}(V)\overline{\circ} \widetilde{\sigma}_B}\arrow{e,t}{\widetilde{\sigma}_V}\node{\widetilde{\sigma}_A\overline{\circ} G^{AB}(V)}
\end{diagram}
$$commutes.

\subsection{The Monoidal bicategory Bieq($\C$) of a tensor category}\label{Bieq(C)}

Given a pair of bicategories $\bb$ and $\bb'$, we can define the ``functor bicategory"  $[\bb,\bb']$, whose 0-cells are pseudofunctors  $\bb\to \bb'$, whose 1-cells are pseudonatural equivalence, and whose 2-cells are invertible modifications.
\smallbreak
The bicategory $[\bb,\bb']$ is not usually a 2-category, because composition of 1-cells in $[\bb,\bb']$ involves composition of 1-cells in $\bb'$, but in the case that $\bb'$ is a 2-category,  $[\bb,\bb']$ is a 2-category.
\smallbreak
When $\bb=\bb'$, the bicategory $[\bb,\bb]$ will be denoted by $\text{Bieq}(\bb)$, and it has a monoidal structure in the sense of \cite{GPS}. Now, we shall describe the monoidal bicategory $\text{Bieq}(\C)$ associated to a tensor category $(\C,\otimes,I)$.
\smallbreak
The tensor product $\overline{\otimes}$ of monoidal endofunctors is defined by the composition of monoidal functors. If $(\theta',\theta'_{(-)}): K\to K'$, $(\theta,\theta_{(-)}):H\to H'$ are pseudonatural transformations, the tensor product $(\theta',\theta'_{(-)})\overline{\otimes} (\theta,\theta_{(-)}): KH\to K'H'$ is defined as $(\theta',\theta'_{(-)})\overline{\otimes} (\theta,\theta_{(-)}):=(K(\theta)\otimes \theta',\theta'\overline{\otimes}\theta)$, where $\theta'\overline{\otimes}\theta$ is given by the commutativity of the following diagram
$$
\begin{diagram}
\node{ KH(V)\otimes K(\theta)\otimes \theta'} \arrow[2]{e,t}{(\theta\overline{\otimes}\theta')_V}\node[2]{K(\theta)\otimes \theta' \otimes K'H'(V)}\\ \node{K(H(V)\otimes \theta)\otimes \theta'}\arrow{n,l}{\psi^K_{H(V),\theta}\otimes\id_{\theta'}}  \arrow{e,t}{K(\theta_V)\otimes\id_{\theta'}}\node{K(\theta\otimes H'(V))\otimes \theta' } \arrow{e,t}{\psi^K_{\theta,H'(V)}\otimes\id_{\theta'}} \node{K(\theta)\otimes K(H'(V))\otimes \theta'} \arrow{n,r}{\id_{K(\theta)}\otimes \theta'_{H'(V)}}
\end{diagram}
$$The tensor product of modifications $g$ and $f$ is defined as $g\overline{\otimes} f := K(f)\otimes g$.
\smallbreak
If  $\chi_1=(\theta,\theta_{(-)}):F\to F'$ and $\chi_2=(\theta',\theta'_{(-)}): H\to H'$ are pseudonatural transformations, where $F,F',H,H':\C\to \C$ are monoidal functors, then there is a natural isomorphism
$$
\begin{diagram}
\node{}\node{F'\circ H}\arrow{se,l}{\id_{F'}\overline{\otimes}\chi_2}\node{}\\
\node{F\circ H}\arrow{ne,l}{\chi_1\overline{\otimes} \id_H}\arrow{se,r}{\id_F\overline{\otimes}\chi_2}\node{\Downarrow c_{\chi_1,\chi_2}}\node{F'\circ H'}\\
\node{}\node{F\circ H'}\arrow{ne,r}{\chi_1\overline{\otimes}\id_H'}
\end{diagram}
$$given by $$c_{\chi_1,\chi_2}:=\theta_{\theta'}^{-1}: \theta\otimes F'(\theta')\to F(\theta')\otimes \theta,$$this natural isomorphism is called the comparison constraint.
\smallbreak
The constraint of associativity $a_{f,g,h}:(f\otimes g)\otimes h\to f\otimes (g\otimes h)$ of the tensor product of pseudonatural transformations $f:K\to K'$, $g:H\to H'$, $h:G\to G'$ is given by the modification $$\psi^K_{H(h),g}\otimes \id_f:KH(h)\otimes K(g)\otimes f\to K(H(h)\otimes g)\otimes f,$$ and it is easy to see that $a$ satisfies the pentagonal identity.
\begin{obs}
The data (TD6), (TD7), and (TD8) of \cite{GPS}, in the monoidal  bicategory $\text{Bieq}(\C)$ are trivial, since we only consider  monoidal functor $(F,\psi):\C\to \C$ such that $F(I)=I$ and $\psi_{V,I}=\psi_{I,V}=\id_V$, for all $V\in \C$.
\smallbreak
The category $\text{Bieq}(\C)(\id_\C,\id_\C)$ is exactly the center of $\C$, \textit{i.e.}, the braided monoidal category $\mathcal{Z}(\C)$, see \cite[pag. 330]{Kassel}.
\end{obs}
\smallbreak
\subsection{Categorical-groups}\label{categorical-groups}

A categorical-group $\mathcal{G}$ is a monoidal category where every object, and every arrow is invertible, \emph{i.e.} for all $X\in \Obj(\mathcal{G})$ there is $X^*\in \Obj(\mathcal{G})$, such that $X\otimes X^*\cong X^*\otimes X\cong I$. We refer the reader to \cite{Baez} for a detailed exposition on the subject.
\smallbreak

A trivial example of a categorical-group is the discrete categorical-group $\underline{G}$, associated to a group $G$. The objects of $\underline{G}$ are the elements of $G$, the arrows are only the identities, and the tensor product is the multiplication of $G$. A more interesting examples is the following.
\begin{ejem}
Let $G$ be a group,  $A$  a $G$-module, and $\omega\in Z^3(G,A)$ a normalized
3-cocycle. We shall define the category  $\C(G,A,\omega)$ by:
\begin{enumerate}
  \item $\Obj(\C(G,A, \omega))=G$,
  \item $\Hom_{\C(G,A, \omega)}(g,h)=\left\{
                     \begin{array}{ll}
                       A, & \hbox{if $g=h$} \\
                       \emptyset, & \hbox{if $g\neq h$.}
                     \end{array}
                   \right.$ \end{enumerate}
We define a  monoidal structure in  $\C(G,A,\omega)$
as follows:

Let $g\in \End(a)$ and $h\in \End(b)$, $a, b \in A$, $g, h \in
G$. Then, $a\otimes b= a+ \  ^gb$ and $g\otimes h = gh$.
We define the associator as $\Phi_{g,h,k}= \omega(g,h,k)$.

The 3-cocycle condition is equivalent to the pentagon axiom, and the condition of normality implies that $e$ is the unit object for this category.
\end{ejem}

\smallbreak
Complete invariants of a categorical-group $\mathcal{G}$ with respect to monoidal equivalences
are
$$\pi_0(\mathcal{G}), \pi_1(\mathcal{G}), \phi(\mathcal{G}),$$ where $\pi_0(\mathcal{G})$ is the group of isomorphism classes of objects, $\pi_1(\mathcal{G})$ is the abelian group of automorphisms of the unit object (the group $\pi_1(\mathcal{G})$ is a $\pi_0(\mathcal{G})$-module in the natural way), and $\phi(\mathcal{G})\in H^3(\pi_0(\mathcal{G}),\pi_1(\mathcal{G}))$ is a third cohomology class given by the associator, see \cite[Subsection 8.3]{Baez} for details on how to obtain $\phi(\mathcal{G})$.
\smallbreak
If $\mathcal{G}$ is a categorical group by \cite[Theorem 43]{Baez} there is an equivalence of monoidal categories between $\mathcal{G}$ and $\C(\pi_0(\mathcal{G}),\pi_1(\mathcal{G}),\phi)$, where $\phi$ is a 3-cocycle in the class  $\phi(\mathcal{G})$.
\smallbreak
Also, it is easy to see that there is a bijective correspondence between monoidal functors  $$F : \C(G,A,\omega) \to \C(G',A',\omega')$$ and triples
$(\pi_0(F),\pi_1(F),\theta(F))$ that consist of:
\begin{itemize}
\item a group morphism $\pi_0(F) : G \to G'$,
\item a $G$-module morphism $\pi_1(F) \colon A \to A'$,
\item a normalized 2-cochain  $k(F) \colon G^2 \to A'$,
such that $dk(F) = \pi_1(F) \omega - \omega' \pi_0(F)^3$.
\end{itemize}
For monoidal functors  $F,F' : \C(G,A,\omega) \to \C(G',A',\omega')$, there is bijective correspondence between monoidal natural isomorphisms $\theta \colon F \to F'$ and normalized 1-cochains
$p(\theta) \colon G \to A'$, where $dp(\theta) = k(F) - k(F')$.

The next result follows from the last discussion or from \cite[Theorem 43]{Baez}.
\begin{prop}\label{prop cat-groups}
Let $\mathcal{G}$ be a categorical group and let $f:G\to \pi_0(\mathcal{G})$ be a morphism of groups. Then there is a monoidal functor $F:\underline{G}\to \mathcal{G}$, such that $f=\pi_0(F)$ if and only if the cohomology class of $\phi f^3$ is zero, where $\phi$ is a 3-cocycle in the class $\phi(\mathcal{G})$.

If  $\phi f^3$ is zero, the classes  of equivalence of monoidal functors $F:\underline{G}\to \mathcal{G}$ such that $\pi_0(F)=f$ are in one to one correspondence with $H^2(G,\pi_1(\mathcal(\mathcal{G})))$.
\end{prop}
\qed
\subsection{Crossed product tensor categories}\label{defincion graded an crossed}

Let $G$ be a group and let $\C$ be a tensor category. We
shall say that  $\C$ is  $G$-graded, if there is a decomposition
$$\C=\oplus_{x\in G}\C_x$$ into a direct sum of full
abelian subcategories, such that  for all $\sigma, x\in G$, the
bifunctor  $\otimes$ maps $\C_\sigma\times \C_x$ to $\C_{\sigma
x}$, see \cite{ENO}.

\begin{defin}
Let $\C$ be a tensor category graded over a group $G$. We shall say that $\C=\bigoplus_{\sigma \in G}\C_\sigma$ is a crossed product tensor category over $G$, if every $\C_\sigma$ has a multiplicatively invertible object.
\end{defin}

Given a group $G$, we define the 2-category of $G$-crossed product tensor categories. The 0-cells are crossed product tensor categories over $G$, 1-cells are graded monoidal functors, \textit{i.e.}, monoidal functors $F: \C\to \D$ such that $F$ maps $\C_\sigma$ to $\D_\sigma$ for all $\sigma \in G$, and 2-cells are  monoidal natural transformations between the graded monoidal functors. The composition of 1-cells and 2-cells is the obvious.

\begin{obs}
The existence of some extra  properties of a  crossed product tensor category $\C$, can be deduced  from the tensor subcategory  $\C_e$. For example $\C$ is semisimple or rigid if and only if $\C_e$ is  semisimple or rigid. However, if $\C_e$ is a braided tensor category, not necessary $\C$ is braided, see Section \ref{seccion 6}.
\smallbreak
A crossed product tensor category $\C$ is a fusion category \cite{ENO} or finite tensor category \cite{finite-categories}, if and only if $G$ is finite, and $\C_e$ is a fusion category or a finite tensor category, respectively.
\end{obs}

\section{Outer $G$-actions over tensor categories}\label{section action y correspondecia}

\subsection{Incoherent outer $G$-actions}

Let $\C$ be a  tensor category.   We define the  categorical-group $\underline{2Out_\otimes(\C)}$, where objects are monoidal autoequivalences  of $\C$, and arrows are  equivalence classes of invertible pseudonatural isomorphisms up to invertible modifications. The composition of arrows in $\underline{2Out_\otimes(\C)}$ is the equivalence class of pseudonatural isomorphisms composition,  and the tensor product is the composition of monoidal functors and pseudonatural transformations.

\begin{defin}
Let $G$ be a group and let $\C$ be a monoidal category. An incoherent outer $G$-action over $\C$, is an \emph{op-monoidal functor} $*:\underline{G}\to \underline{2Out_\otimes(\C)}$. Two incoherent outer $G$-actions are equivalent if the associated monoidal functors are monoidally equivalent.
\end{defin}

We shall analyze the incoherent outer $G$-action using the Subsection \ref{categorical-groups}. Complete invariants for the categorical group $\underline{2Out_\otimes(\C)}$ are $\pi_0(\underline{2Out_\otimes(\C)})$ the equivalences classes of monoidal functor under invertible modification,  $\pi_1(\underline{2Out_\otimes(\C)})=Inv(\mathcal{Z}(\C))$ the abelian group of isomorphisms classes of invertible objects of the center of $\C$, and a third cohomology class $\phi(\underline{2Out_\otimes(\C)}) \in H^3(\pi_0(\underline{2Out_\otimes(\C)}),Inv(\mathcal{Z}(\C)))$.
\smallbreak

Every incoherent outer $G$-action over a tensor category induces a group morphism $f : G \to \pi_0(\underline{2Out_\otimes(\C)})$. We shall say that a group morphism $f : G \to \pi_0(\underline{2Out_\otimes(\C)})$ is realizable if there is some incoherent outer $G$-action such that the induced group morphism coincides with $f$.

\begin{prop}
Let $G$ be a group and let $f:G\to \pi_0(\underline{Out_\otimes(\C)})$ be a group morphism. Then there is an incoherent outer $G$-action over $\C$ that realize the morphism $f$ if and only if the cohomology class of $\phi f^3$ is zero, where $\phi$ is some 3-cocycle in the class of $\phi(\underline{2Out_\otimes(\C)})$.

If  $\phi f^3$ is zero, the classes of equivalence of monoidal functors $F:\underline{G}\to \underline{2Out_\otimes(\C)}$ such that $\pi_0(F)=f$ are in one to one correspondence with $H^2(G, Inv(\mathcal{Z}(\C)))$.
\end{prop}
\begin{proof}
See Proposition \ref{prop cat-groups}.
\end{proof}

\subsection{Coherent outer $G$-actions}

Let $\C$ be a monoidal category and let $F: G\to \underline{2Out_\otimes(\C)}$ be an incoherent outer $G$-action. We define a crossed system  associated to $F$ as the following data

\begin{itemize}
  \item monoidal functors $(\sigma_*,\psi^{\sigma_*}): \C\to \C$ for all $\sigma \in G$,
  \item pseudonatural isomorphisms $(U_{\sigma,\tau},\chi_{\sigma,\tau}): \sigma_*\circ \tau_*\to (\sigma\tau)_*$ for all $\sigma,\tau  \in G$,
  \item invertible modifications $\omega_{\sigma,\tau,\rho}: \chi_{\sigma,\tau\rho}\overline{\circ} (\id_{\sigma_*} \overline {\otimes}\chi_{\tau,\rho}) \to \chi_{\sigma\tau,\rho}\overline{\circ} (\chi_{\sigma,\tau}\overline{\otimes}\id_{\rho_*})$ for all $\sigma,\tau,\rho \in G$.
\end{itemize}
that realize the incoherent outer $G$-action $F$ (recall that the symbol $\overline{\circ}$ is the composition of 1-cells in the bicategory Bieq$(\C)$).

\begin{obs}
\begin{enumerate}
  \item By abuse of notation, we write $\chi_{\sigma,\tau}$ instead of pseudonatural transformation $(U_{\sigma,\tau},\chi_{\sigma,\tau})$, when no confusion can arise.

  \item By the definition of $\underline{2Out_\otimes(\C)}$, there are several crossed systems that realize an incoherent $G$-action.

  \item For every crossed system, without loss of generality, we can and shall assume that

\begin{itemize}
  \item $e_*= \id_\C$\  the  monoidal identity functor,
  \item $\chi_{e,\sigma}=\chi_{\sigma,e}=(I,\id_{\sigma_*})$ the identity pseudonatural isomorphism,
  \item $ \omega_{\sigma,e,\tau}=\id_{\chi_{\sigma,\tau}}$ the identity modification,
\end{itemize}for all $\sigma,\tau \in G$.

\end{enumerate}
\end{obs}

\begin{defin}
Let $\C$ be a tensor category and let $F: \underline{G}\to \underline{2Out_\otimes(\C)}$ be an incoherent outer $G$-action. A coherent outer $G$-action associated to $F$, is a crossed system $(\{ \widetilde{\sigma}\}_{\sigma\in G},\chi,\omega)$ associated to $F$, such that for the pseudonatural isomorphisms $\chi_{\sigma,\tau}$ and the invertible modifications $\omega_{\sigma,\tau,\rho}$, the diagram

\begin{equation}\label{coherencia} \begin{diagram}\node{}\node{\sigma_*(U_{\tau,\rho}U_{\tau\rho,\mu})U_{\sigma,\tau\rho\mu}}\arrow{se,t}{\psi^{\sigma_*}_{U_{\tau,\rho},U_{\tau\rho,\mu}}\id_{U_{\sigma,\tau\rho\mu}}}\node{}\\ \node{[\sigma_*(\tau_*(U_{\rho,\mu})U_{\tau,\rho\mu})]U_{\sigma,\tau\rho\mu}} \arrow{s,r}{\psi^{\sigma_*}_{\tau_*(U_{\rho,\mu}),U_{\tau,\rho\mu}}\id_{U_{\sigma,\tau\rho\mu}}} \arrow{ne,t}{\sigma_*(\omega_{\tau,\rho,\mu})\id_{U_{\sigma,\tau\rho\mu}}} \node{}  \node{[\sigma_*(U_{\sigma,\tau})\sigma_*(U_{\tau\rho,\mu})]U_{\sigma,\tau\rho\mu}} \arrow{s,r}{\alpha_{\sigma_*(U_{\sigma,\tau}),\sigma_*(U_{\tau\rho,\mu}), U_{\sigma,\tau\rho\mu}} }\\ \node{[\sigma_*(\tau_*(U_{\sigma,\mu}))\sigma_*(U_{\tau,\rho\mu})]U_{\sigma,\tau\rho\mu}}\arrow{s,r}{\alpha_{\sigma_*(\tau(U_{\sigma,\tau})),\sigma_*(U_{\tau\rho,\mu}),U_{\sigma,\tau\rho\mu}} }  \node{}  \node{\sigma_*(U_{\sigma,\tau})[\sigma_*(U_{\tau\rho,\mu})U_{\sigma,\tau\rho\mu}]}\arrow{s,r}{\id_{\sigma_*(U_{\sigma,\tau})}\omega_{\sigma,\tau\rho,\mu}}\\ \node{\sigma_*(\tau_*(U_{\sigma,\mu}))[\sigma_*(U_{\tau,\rho\mu})U_{\sigma,\tau\rho\mu}]} \arrow{s,r}{\id_{\sigma_*(\tau_*(U_{\sigma,\tau}))}\omega_{\sigma,\tau,\rho\mu}} \node{}\node{\sigma_*(U_{\tau,\rho})[U_{\sigma,\tau\rho}U_{\sigma\tau\rho,\mu}]} \arrow{s,r}{\alpha_{\sigma_*(U_{\tau,\rho}),U_{\sigma,\tau\rho},U_{\sigma\tau\rho,\mu}}^{-1}}\\ \node{\sigma_*(\tau_*(U_{\rho,\mu}))[U_{\sigma,\tau}U_{\sigma\tau,\rho\mu}]} \arrow{s,r}{\alpha_{\sigma_*(\tau_*(U_{\rho,\mu})),U_{\sigma,\tau},U_{\sigma\tau,\rho\mu}}^{-1}} \node{}\node{[\sigma_*(U_{\tau,\rho})U_{\sigma,\tau\rho}]U_{\sigma\tau\rho,\mu}}\arrow{s,r}{\omega_{\sigma,\tau,\rho}\id_{U_{\sigma\tau\rho,\mu}}}\\ \node{[\sigma_*(\tau_*(U_{\rho,\mu}))U_{\sigma,\tau}]U_{\sigma\tau,\rho\mu}} \arrow{s,r}{\chi_{\sigma,\tau}(U_{\rho,\mu})\id_{U_{\sigma\tau,\rho\mu}}} \node{}  \node{[U_{\sigma,\tau}U_{\sigma\tau,\rho}]U_{\sigma\tau\rho,\mu}}\arrow{s,r}{\alpha_{U_{\sigma,\tau},U_{\sigma\tau,\rho},U_{\sigma\tau\rho,\mu}}^{-1}}\\ \node{[U_{\sigma,\tau}(\sigma\tau)_*(U_{\rho,\mu})]U_{\sigma\tau,\rho\mu}}\arrow{se,t}{\alpha_{U_{\sigma,\tau},(\sigma\tau)_*(U_{\rho,\mu}),U_{\sigma\tau,\rho\mu}}} \node{}\node{U_{\sigma,\tau}[U_{\sigma\tau,\rho}U_{\sigma\tau\rho,\mu}]}\\ \node{}\node{U_{\sigma,\tau}[(\sigma\tau)_*(U_{\rho,\mu})U_{\sigma\tau,\rho\mu}]} \arrow{ne,b}{\id_{U_{\sigma,\tau}}\omega_{\sigma\tau,\rho,\mu}}\node{} \end{diagram} \end{equation} commutes for all $\sigma,\tau,\rho,\mu \in G$ (where tensor symbols among objects and arrows have been omitted as a space-saving measure).

\end{defin}

\begin{obs}\label{3-group discreto}
 For every group $G$, we can associate a discrete 3-category $\underline{\underline{G}}$, where objects are elements of $G$, and

$$\underline{\underline{G}}(g,h)=\left\{
                       \begin{array}{ll}
                         \{*\}, & \hbox{if $g=h$} \\
                         \emptyset, & \hbox{if $g\neq h$.}
                       \end{array}
                     \right.$$
\begin{gather*}\overline{\otimes}: \underline{\underline{G}}\times \underline{\underline{G}}\to \underline{\underline{G}},\\
\underline{\otimes}_{g,h}:\underline{\underline{G}}(\sigma,\sigma)\times \underline{\underline{G}}(\tau,\tau) \to \underline{\underline{G}}(\sigma\tau,\sigma\tau).
\end{gather*}
The definition of a coherence outer $G$-action over $\C$, is equivalent to the definition of a trihomorphism from $\underline{\underline{G}}$ to $\text{Bieq}(\C)$ (see \cite{GPS} for the definition of trihomomorphism).
\end{obs}
Given a crossed system associated to an incoherent outer $G$-action, we can define a monoidal bicategory. In order to describe the monoidal bicategory in a simple way we can suppose, without loss of generality that $\C$ is skeletal, so for every pair $\sigma,\tau \in G$ there is only one  pseudonatural transformation $\chi_{\sigma,\tau}^{-1}:(\sigma\tau)_*\to \sigma_*\tau_*$, such that $\chi_{\sigma,\tau}^{-1} \circ \chi_{\sigma,\tau }= \id_{\sigma_*\tau_*}$, and $\chi_{\sigma,\tau} \circ\chi_{\sigma,\tau}^{-1} = \id_{(\sigma\tau)_*}$ for all $\sigma, \tau\in G$. Let $<G>\subseteq \text{Bieq}(\C)$ be the full sub-bicategory where objects are $\{\sigma_*\}_{\sigma\in G}$. We define a homomorphism of bicategories $\otimes_G: <G> \times <G> \to <G>$ by $\sigma_*\otimes_G \tau_* = (\sigma\tau)_*$, and the commutativity of the diagram
$$
\begin{diagram}
\node{\sigma_*\tau_*} \arrow{s,l}{\chi_{\sigma,\tau}} \arrow{e,t}{f\overline{\otimes} g}\node{\sigma_*\tau_*}\arrow{s,r}{\chi_{\sigma,\tau}}\\
\node{(\sigma\tau)_*}\arrow{e,t}{f\otimes_G g}\node{(\sigma\tau)_*}
\end{diagram}
$$ where $f\in \text{Bieq}(\sigma_*,\sigma_*), g\in \text{Bieq}(\C)(\tau_*,\tau_*)$.

We define a pseudonatural equivalence in the bicategory $[<G>\times <G>\times <G>, <G>]$ by the commutativity of the diagram
$$
\begin{diagram}
\node{(\sigma\tau\rho)_*}\arrow[2]{e,t}{a_{\sigma,\tau,\rho}}\node{} \node{(\sigma\tau\rho)_*}\\
\node{\sigma_*\circ(\tau\rho)_*}\arrow{n,l}{\chi_{\sigma,\tau\rho}} \node{}\node{(\sigma\tau)_*\rho}\arrow{n,r}{\chi_{\sigma\tau,\rho}}\\
\node{}  \node{\sigma_*\tau_*\rho_*} \arrow{nw,r}{\id_\sigma\overline{\otimes}\chi_{\tau,\rho}} \arrow{ne,r}{\chi_{\sigma,\tau}\overline{\otimes}\id_\rho}
\end{diagram}
$$ The diagram \eqref{coherencia} define a  modification $$\pi: a_{\sigma\tau,\rho,\mu}\circ a_{\sigma,\tau,\rho\mu} \to a_{\sigma,\tau,\rho}\otimes_G \id_\mu \circ a_{\sigma,\tau\rho,\mu}\circ \id_{\sigma}\otimes_G a_{\tau,\rho,\mu},$$  in the bicategory $[<G>\times <G>\times <G>\times <G>, <G>]$.

Since for every invertible object $U\in \C$, we can identify $\Aut_\C(U)$ with $\Aut(I)=k^*$,  the associativity constraint $\alpha$, the pseudo-natural isomorphism $\chi$ and the natural isomorphism $\psi^\sigma,$ restricted to the invertible objects $\operatorname{Inv}(\C)$ define the maps 
\begin{align*}
&\alpha: \operatorname{Inv}(\C)\times \operatorname{Inv}(\C)\times \operatorname{Inv}(\C) \to k^*,\\
&\psi^{\sigma}:\operatorname{Inv}(\C)\times \operatorname{Inv}(\C) \to k^*, \  \  (\sigma \in G)\\
&\chi_{\sigma,\tau}: \operatorname{Inv}(\C)\to k^*,\  \  \   (\sigma,\tau \in G),
\end{align*}
that is,

\begin{align*}
&\alpha_{A,B,C}:=\alpha(A,B,C)\id_{ABC}:[AB]C\to A[BC], \  \  (A,B,C\in \operatorname{Inv}(\C))\\
&\psi^{\sigma}_{A,B}:=\psi^{\sigma}(A,B)\id_{\sigma_*(AB)}: \sigma_*(AB)\to \sigma_*(A)\sigma_*(B) \  \  (\sigma \in G, A,B\in \operatorname{Inv}(\C))\\
&\chi_{\sigma,\tau}(A):= \sigma_*(\tau_*(A))U_{\sigma,\tau}\to U_{\sigma,\tau}(\sigma\tau)_*(A),\  \  \   (\sigma,\tau \in G, A\in \operatorname{Inv}(\C)).
\end{align*}
 We can define a map  $$\pi:G\times G\times G\times G\to k^*,$$ as the error of the commutativity of diagram \eqref{coherencia}:

\begin{align}
\pi(\sigma,\tau,\rho,\mu)=&\theta_{\sigma,\tau}(U_{\rho,\mu})\\
&c((\sigma\tau)_*(U_{\rho,\mu}),U_{\sigma,\tau})\notag\\
&\psi^{\sigma_*}(\tau_*(U_{\rho,\mu}),U_{\tau,\rho\mu})\notag\\ &\psi^{\sigma_*}(U_{\tau,\rho},U_{\tau\rho,\mu})^{-1}\notag \\
&\alpha(\sigma_*(\tau_*(U_{\sigma,\tau}) ),\sigma_*(U_{\tau\rho,\mu}),U_{\sigma,\tau\rho\mu})\notag\\
&\alpha(\sigma_*(\tau_*(U_{\rho,\mu})),U_{\sigma,\tau},U(\sigma\tau,\rho\mu))^{-1}\notag\\
&\alpha(U_{\sigma,\tau},(\sigma\tau)_*(U_{\rho,\mu}),U_{\sigma,\tau,\rho\mu})\notag\\
&\alpha(U_{\sigma,\tau},U_{\sigma\tau,\rho},U_{\sigma\tau\rho,\mu})^{-1}\notag\\
&\alpha(\sigma_*(U_{\tau,\rho}),U_{\sigma,\tau\rho},U_{\sigma\tau\rho,\mu})\notag\\
&\alpha(\sigma_*(U_{\sigma,\tau}),\sigma_*(U_{\tau\rho,\mu}),U_{\sigma,\tau\rho\mu})^{-1}.\notag
\end{align} It is  straightforward  (but tedious and long) to see that $\pi$ is a 4-cocycle, see \cite[Subsection 8.4]{ENO3}. It is also possible to see the 4-cocycle condition directly for the nonabelian 4-cocycle condition \cite[(TA1)]{GPS} in the monoidal bicategory $<G>$.  It is also straightforward  to see that if the chosen of the crossed system is changed (by an equivalent one), the 4-cocycle $\pi$ only change for a 4-coboundary, so an incoherent outer $G$-action defines a fourth cohomology class.

\begin{prop}
Let $\C$ be a tensor category. An incoherent outer $G$-action over $\C$ is coherent if and only if the  associated fourth cohomology class is trivial.
\end{prop}
\begin{proof}
If an  outer $G$-action is coherent, the diagram \eqref{coherencia} commutes, so the map $\pi$ is trivial. Conversely, if there is a 3-coboundary  $\lambda: G\times G\times G\to G$, such that $\delta(\lambda)= \pi$, then the  modification defined by the map $\lambda^{-1}\omega$ defines a coherent outer $G$-action.
 \end{proof}

\subsection{Crossed product tensor category associated to a coherent outer $G$-action}\label{from action to crossed}

If a group $G$ acts over a monoidal category $\C$, we shall define a $G$-crossed product tensor category associated to this action, denoted as $\C\rtimes G$. We set $\C\rtimes G=\bigoplus_{\sigma}\C_\sigma$ as an abelian category, where $\C_\sigma =\C$. We shall denote by $[V,\sigma]$ the object $V\in \C_{\sigma}$, and a morphism from $\bigoplus_{\sigma\in G}[V_\sigma,\sigma]$ to $\bigoplus_{\sigma\in G}[W_\sigma,\sigma]$ is expressed as $\bigoplus_{\sigma\in G}[f_\sigma,\sigma]$ with $f_\sigma: V_\sigma\to W_\sigma$ a morphism in $\C$.
\smallbreak
The tensor product  $\cdot:\C\rtimes G\times \C\rtimes G\to \C\rtimes G$ is defined by
\begin{align*}
    [V,\sigma]\cdot [W,\tau] &:= [V\otimes\sigma_*(W)\otimes U_{\sigma,\tau},\sigma\tau] \  \  \text{ for objects, and}\\
    [f,\sigma]\cdot [g,\tau] &:= [f\otimes \sigma_*(g)\otimes \id_{U_{\sigma,\tau}},\sigma\tau]\  \  \text{ for morphisms.}
\end{align*}

It is easy to see that the unit object is $(I,e)$. The associativity is given by

$$
\begin{diagram}
\node{[V,\sigma]\cdot([W,\tau]\cdot[Z,\rho])} \arrow{e,t,=}{}\arrow{s,l}{\alpha_{[V,\sigma],[W,\tau],[Z,\rho]}} \node{[V\otimes\sigma_*(W\otimes\tau_*(Z)\otimes U_{\tau,\rho})\otimes U_{\sigma,\tau\rho},\sigma\tau\rho]} \arrow{s,r}{[\alpha^{\C\rtimes G}(V,\sigma,W,\tau,Z,\rho),\sigma\tau\rho]}\\
\node{([V,\sigma]\cdot[W,\tau])\cdot[Z,\rho]}\arrow{e,=}{} \node{[V\otimes\sigma_*(W)\otimes U_{\sigma,\tau}\otimes(\sigma\tau)_*(Z)\otimes U_{\sigma\tau,\rho},\sigma\tau\rho]}
\end{diagram}
$$

where $\alpha^{\C\rtimes G}(V,\sigma,W,\tau,Z,\rho)$ is the composition
$$
\begin{diagram}
\node{V\otimes\sigma_*(W\otimes\tau_*(Z)\otimes U_{\tau,\rho})\otimes U_{\sigma,\tau\rho}}\arrow{s,l} {\id_V\otimes\psi^{\sigma_*}_{W,\tau_*(Z)\otimes U_{\tau,\sigma}}\otimes \id_{U_{\sigma,\tau\rho}}}\\
\node{V\otimes\sigma_*(W)\otimes\sigma_*(\tau_*(Z)\otimes U_{\tau,\rho})\otimes U_{\sigma,\tau\rho}}  \arrow{s,l}{\id_{V\otimes\sigma_*(W)}\otimes \psi^{\sigma_*}_{\tau_*(Z),U_{\tau,\rho}}\otimes \id_{U_{\sigma,\tau\rho}}}\\
\node{V\otimes\sigma_*(W)\otimes\sigma_*(\tau_*(Z))\otimes \sigma_*(U_{\tau,\rho})\otimes U_{\sigma,\tau\rho}} \arrow{s,l}{\id_{{V\otimes\sigma_*(W)\otimes\sigma_*(\tau_*(Z))}}\otimes \omega_{\sigma,\tau,\rho}}\\
\node{V\otimes\sigma_*(W)\otimes\sigma_*(\tau_*(Z))\otimes U_{\sigma,\tau}\otimes U_{\sigma\tau,\rho}} \arrow{s,l}{\id_{V\otimes \sigma_*(W)}\otimes\chi_{\sigma,\tau}\otimes\id_{U_{\sigma\tau,\rho}}}\\
\node{V\otimes\sigma_*(W)\otimes U_{\sigma,\tau}\otimes (\sigma\tau)_*(Z)\otimes U_{\sigma\tau,\rho}}
\end{diagram}
$$
The associativity constraint have been omitted as a space-saving measure.
As we shall see, the coherence  condition over an outer $G$-action, is exactly the pentagonal identity for $\C\rtimes G$.

\subsubsection{Pentagonal identity for $\C\rtimes G$}

For a category $\D$ with a bifunctor $\otimes:\D\times\D\to \D$ and natural isomorphisms $\alpha_{A,B,C}: A\otimes(B\otimes C)\to (A\otimes B)\otimes C$, we shall denote by $P(A,B,C,D)$ the following pentagonal diagram
$$
\begin{diagram}
\node{}\node{A\otimes(B\otimes(C\otimes D))} \arrow{sw,l}{\alpha_{A,B,C\otimes D}}\arrow{se,l}{\id\otimes \alpha_{B,C,D}}\node{}\\
\node{(A\otimes B)\otimes (C\otimes D)}\arrow{s,l}{\alpha_{A\otimes B,C,D}} \node{}\node{A\otimes((B\otimes C)\otimes D)} \arrow{s,r}{\alpha_{A,B\otimes C,D}}\\
\node{((A\otimes B)\otimes C)\otimes D} \node{} \node{(A\otimes (B\otimes C))\otimes D} \arrow[2]{w,t}{\alpha_{A,B,C}\otimes \id}
\end{diagram}
$$

\begin{obs}\label{obs notacion}
From now on, we shall denote $[V]:= [V,e]$ and $[\sigma]:=
[I,\sigma]$, for all $V\in \C, \sigma \in G$. Analogously, $[f]:=
[f,e]: [V]\to [W]$  for all arrow $f:V\to W$ in $\C$. Note  that
$[V]\cdot[\sigma]= [V,\sigma]$ and
$[\sigma]\cdot[V]=[\sigma_*(V),\sigma]$ for all $V\in \C, \sigma \in
G$.
\end{obs}

In order to prove the coherence of $\C\rtimes G$, is sufficient to see the pentagonal identity for the $[V], [\sigma]$, $V\in \C, \sigma \in G$, since every object in $\C\rtimes G$ is a direct sum of tensor products of $[V], [\sigma]$.
\smallbreak
First, see the next equality

\begin{equation}\label{igualdad asociatividad}
\alpha_{[V,\sigma],[W,\tau],[Z,\rho]}=\id_{[V]}\cdot\alpha_{[I,\sigma],[W,\tau],[Z,\rho]},
\end{equation} so $\alpha_{[V,e],[W,\tau],[Z,\rho]}=\id$. The eight pentagonal identities
\begin{multline*}
P([V],[\sigma],[\tau],[\rho]), P([V],[\sigma],[W],[\tau]), P([V],[W],[\sigma],[Z]), P([V],[W],[Z],[\sigma]),\\  P([V],[W],[\sigma],[\tau]), P([V],[\sigma],[W],[Z]),   P([V],[\sigma],[\tau],[W]), P([V],[W],[Z],[U])
\end{multline*}
follow  from \eqref{igualdad asociatividad}.
\smallbreak
The pentagon $P([\sigma],[V],[W],[\tau])$ commutes because $$\alpha_{[\sigma],[W],[Z,\rho]}=\alpha_{[\sigma],[W],[Z]}\cdot\id_{[\rho]}.$$
\smallbreak
The pentagons $P([\sigma],[V],[\tau],[W]), P([\sigma],[\tau],[\rho])$, $P([\sigma],[\tau],[V],[\rho])$ commute by the definition of $\alpha$.
\smallbreak
The following table explains the commutativity of the other pentagons

\begin{tabular}{|c|c|}
  \hline
  Pentagons &  Pentagonal identity equivalence\\
  \hline
  $P([\sigma],[V],[W],[Z])$ & $(\sigma_*,\psi^{\sigma_*})$ is a monoidal functor\\
  \hline
  $P([\sigma],[\tau],[V],[W])$ & $(\chi_{\sigma,\tau},U_{\sigma,\tau})$ is a pseudonatural equivalence\\
  \hline
  $P([\sigma],[\tau],[\rho],[V])$& $\omega_{\sigma,\tau,\rho}$ is a modification\\
  \hline
  $P([\sigma],[\tau],[\rho],[\gamma])$ & commutativity of the diagram \eqref{coherencia}.\\
  \hline
\end{tabular}

\subsection{The coherent outer $G$-action associated to a $G$-crossed product tensor category}

Let $\C$ be a $G$-crossed product tensor category. In order to show more clearly the associated coherent outer $G$-action, we shall make some reductions. Let we choose a family $\{N_\sigma\}_{\sigma\in G}$ of homogeneous invertible objects, where $N_e=I$. The family $\{N_\sigma\}_{\sigma\in G}$ defines the equivalences of categories
\begin{align*}
    N_\sigma\otimes (-): \C_e &\to \C_\sigma\\
     V &\mapsto V\otimes N_\sigma\\
     f&\mapsto f\otimes \id_{N_\sigma}.
\end{align*}
Using these equivalences, we have an equivalence of categories $$\C=\bigoplus_{\sigma\in G}\C_\sigma\to \bigoplus_{\sigma\in G}\C_e^\sigma,$$ where $\C_e^\sigma=\C_e$, for all $\sigma\in G$.
\smallbreak
Now, we  can transport the monoidal structure of $\C$ to $\bigoplus_{\sigma\in G}\C_e$. Then, without loss of generality we can suppose that the graded tensor category $\C=\bigoplus_{\sigma\in G}\C_\sigma$ has the following properties:

\begin{itemize}
  \item $\C_\sigma=\C_e$ for all $\sigma \in G$ (so we can and will use the same notations of the Remark \ref{obs notacion}),
  \item the objects $[\sigma]\in \C_\sigma$ are invertible for all $\sigma\in G$,
  \item $[V]\cdot[W,\sigma]=[V\otimes W,\sigma]$, for all $V, W \in \C_e$, $\sigma\in G$.
\end{itemize}

For each pair $\sigma, \tau \in G$, we have that $[\sigma]\cdot [\tau] \in \C_{\sigma\tau}$, so there is an unique invertible object $U_{\sigma\tau}\in \C_e$, such that $[\sigma]\cdot [\tau]=[U_{\sigma,\tau},\sigma\tau]$. Analogously, the  objects $[\sigma]$ define  functors $\sigma_*:\C_e\to \C_e, V\mapsto \sigma_*(V)$ by the rule $[\sigma]\cdot [V]=[\sigma_*(V),\sigma]$ for all $V\in \C_e$, and $\id_{[\sigma]}\cdot [f] =[\sigma_*(f),\sigma]$, for all arrow $f$ in $\C_e$.

\begin{lem}\label{lema producto}
If  the category $(\C_e,\otimes, I)$ is skeletal, then
$$[V,\sigma]\cdot[W,\tau]=[V\otimes\sigma_*(W)\otimes U_{\sigma,\tau},\sigma\tau]$$
for all $V, W\in \C_e, \   \sigma, \tau \in G$.
\end{lem}
\begin{proof}
Since $\C_e$ is skeletal, the category $\C=\bigoplus_{\sigma}\C_\sigma$  is skeletal. Then we do not need to parenthesize tensor
products for objects in $\C$. Also, recall that $[\sigma]\cdot[V]=[\sigma_*(V),\sigma]=[\sigma_*(V)]\cdot[\sigma]$, for all $V\in \C_e$, $\sigma \in G$.

Hence,
\begin{align*}
    [V,\sigma]\cdot[W,\tau] &= [V]\cdot[\sigma]\cdot[W]\cdot[\tau]\\
    &=[V]\cdot[\sigma_*(W)]\cdot[\sigma]\cdot[\tau] \\
    &=[W\otimes \sigma_*(W)]\cdot[U_{\sigma,\tau},\sigma\tau]\\
    &=[W\otimes \sigma_*(W)]\cdot[U_{\sigma,\tau}]\cdot[\sigma\tau]\\
    &=[W\otimes \sigma_*(W)\otimes U_{\sigma,\tau}]\cdot[\sigma\tau]\\
    &=[W\otimes \sigma_*(W)\otimes U_{\sigma,\tau},\sigma\tau]\\
\end{align*}for all $V, W \in \C_e$, $\sigma \in G$.
\end{proof}
Under this reduction and using the Lemma \ref{lema producto} we can describe the coherent outer $G$-action as the reciprocal construction of the Subsection \ref{from action to crossed}. Suppose that $\C_e$ is skeletal, then the data that define the coherent outer $G$-action associated to $\C$ are the following:

\begin{itemize}
  \item monoidal equivalences: $(\sigma_*,\psi^{\sigma_*}): \C_e\to \C_e$, where $$[\psi^{\sigma_*}_{W,Z},\sigma]:=\alpha_{[\sigma],[W],[Z]}:[\sigma_*(W\otimes Z),\sigma]\to [\sigma_*(W)\otimes \sigma_*(Z),\sigma],$$
  \item pseudonatural transformations: $ (U_{\sigma,\tau},\chi_{\sigma,\tau}):\sigma_*\circ \tau_*\to (\sigma\tau)_*$, where $$[\chi_{\sigma,\tau}(Z),\sigma\tau]:=\alpha_{[\sigma],[\tau],[Z]}: [\sigma_*(\tau_*(Z))\otimes U_{\sigma,\tau},\sigma\tau] \to [U_{\sigma,\tau}\otimes(\sigma\tau)_*(Z),\sigma\tau],$$
  \item modifications $\omega_{\sigma,\tau\rho}:  \chi_{\sigma,\tau\rho}\overline{\circ} (\id_{\sigma_*} \overline {\otimes}\chi_{\tau,\rho}) \to \chi_{\sigma\tau,\rho}\overline{\circ} (\chi_{\sigma,\tau}\overline{\otimes}\id_{\rho_*})$, where
$$[\omega_{\sigma,\tau,\rho},\sigma\tau\rho]:= \alpha_{[\sigma],[\tau],[\rho]}:[\sigma_*(U_{\tau,\rho})\otimes U_{\sigma,\tau\rho},\sigma\tau\rho] \to [U_{\sigma,\tau}\otimes U_{\sigma\tau,\rho},\sigma\tau\rho].$$
\end{itemize}

\section{$G$-crossed product tensor category in terms of  coherent outer $G$-actions}\label{section classif}

In this section we shall define the 2-category of coherent outer $G$-actions, and  we shall see that the 2-category of crossed product tensor category over a fixed group $G$, is equivalent to the 2-category of all coherent outer $G$-actions.
\smallbreak
The 0-cells of the 2-category of coherent outer $G$-actions are coherent outer $G$-action over a tensor category.
\smallbreak
Let $(\{ \widetilde{\sigma}\}_{\sigma\in G},\chi,\omega)$ and $(\{\widehat{\sigma}\}_{\sigma \in G},\chi',\omega')$ be coherent outer $G$-actions over tensor categories $\C$ and $\D$, respectively. An arrow from $(\{ \widetilde{\sigma}\}_{\sigma\in G},\chi,\omega)$ to $(\{\widehat{\sigma}\}_{\sigma \in G},\chi',\omega')$, is a triple $(H,\theta,\Pi)$, where $(H,\psi^H):\C\to \D$ is an monoidal functor,  $(\theta^\sigma,\theta_\sigma): \widehat{\sigma}\circ H\to H\circ \widetilde{\sigma}$ is a pseudonatural equivalence for each $\sigma\in G$, and $\Pi$ is  a modification

$$
\begin{diagram}
\node{}\node{\widehat{\sigma} H\widetilde{\tau}}\arrow{se,t}{\theta_\sigma\overline{\otimes}id_{\widetilde{\tau}} } \node{}\\
\node{\widehat{\sigma}\widehat{\tau}H}\arrow{ne,t}{\id_{\widehat{\sigma}}\overline{\otimes}\theta_\tau} \arrow{s,t}{\chi_{\sigma,\tau}'\overline{\otimes}\id_H}\node{\Uparrow\Pi_{\sigma,\tau}} \node{H \widetilde{\sigma}\widetilde{\tau}}\arrow{s,r}{\id_H\overline{\otimes} \chi_{\sigma,\tau}}\\
\node{\widehat{\sigma\tau}H}\arrow[2]{e,t}{\theta_{\sigma\tau}}\node{} \node{H\widetilde{\sigma\tau}}
\end{diagram}
$$
such that: $(\theta^e, \theta_e)=(I,\id)$, $\Pi_{\sigma,e}=\Pi_{e,\sigma}=\id_{\theta_\sigma}$ for all $\sigma \in G$, and the diagram

$$
\begin{diagram}
\node{\theta_\sigma\widehat{\sigma}(\theta_\tau\widehat{\tau}(\theta_\rho)U_{\tau,\rho}')U_{\sigma,\tau\rho}'} \arrow{s,l}{\id_{\theta_\sigma}\otimes\widehat{\sigma}(\Pi_{\tau,\rho})\otimes\id_{U_{\sigma,\tau\rho}'}} \arrow{e,t}{\alpha^{\D\rtimes G}(\theta_\sigma,\sigma,\theta_\tau,\tau,\theta_\rho,\rho)} \node{\theta_\sigma\widehat{\sigma}(\theta_\tau)U_{\sigma,\tau}'\widehat{\sigma\tau}(\theta_\rho)U_{\sigma\tau,\rho}'} \arrow{s,r}{\Pi_{\sigma,\tau}\otimes\id_{\widehat{\sigma\tau}(\theta_\rho)U_{\sigma\tau,\rho}'}}\\
\node{\theta_\sigma\widehat{\sigma}(H(U_{\tau,\rho})\theta_{\tau\rho})U_{\sigma,\tau\rho}'} \arrow{s,l}{\id_{\theta_\sigma}\otimes\psi^{\widehat{\sigma}}_{H(U_{\tau,\rho}),\theta_{\tau\rho}}\otimes\id_{U_{\sigma,\tau\rho}'}} \node{H(U_{\sigma,\tau})\theta_{\sigma\tau}\widehat{\sigma\tau}(\theta_\rho)U_{\sigma\tau,\rho}'} \arrow{s,r}{\id_{H(U_{\sigma,\tau})}\otimes\Pi_{\sigma\tau,\rho}}\\
\node{\theta_\sigma\widehat{\sigma}(H(U_{\tau,\rho}))\widehat{\sigma}(\theta_{\tau\rho})U_{\sigma,\tau\rho}'} \arrow{s,l}{\theta^\sigma_{U_{\tau,\rho}}\otimes\id_{\widehat{\sigma}(\theta_{\tau\rho})U_{\sigma,\tau\rho}'}} \node{H(U_{\sigma,\tau})H(U_{\sigma\tau,\rho})\theta_{\sigma\tau\rho}} \arrow{s,r}{\psi^H_{U_{\sigma,\tau},U_{\sigma\tau,\rho}}\otimes\id_{\theta_{\sigma\tau\rho}}}\\
\node{H(\widetilde{\sigma}(U_{\tau,\rho}))\theta_\sigma\widehat{\sigma}(\theta_{\tau\rho})U_{\sigma,\tau\rho}'} \arrow{s,l}{\id_{H(\widetilde{\sigma}(U_{\tau,\rho}))}\otimes\Pi_{\sigma,\tau\rho}} \node{H(U_{\sigma,\tau}U_{\sigma\tau,\rho})\theta_{\sigma\tau\rho}}\\
\node{H(\widetilde{\sigma}(U_{\tau,\rho}))\otimes H(U_{\sigma,\tau\rho})\theta_{\sigma\tau\rho}} \arrow{e,t}{\psi^H_{\widetilde{\sigma}(U_{\tau,\rho}),U_{\sigma,\tau\rho}}} \node{H(\sigma_*(U_{\tau,\rho})U_{\sigma,\tau\rho})\theta_{\sigma\tau\rho}} \arrow{n,r}{H(\omega_{\sigma,\tau,\rho})}
\end{diagram}
$$commutes for all $\sigma,\tau, \rho \in G$ (where tensor symbols among objects  have been omitted as a space-saving measure).
\smallbreak
A 2-cell from $(H,\theta,\Pi)$ to $(\widetilde{H},\widetilde{\theta},\widetilde{\Pi})$ consist of the data $\{m_\sigma, m\}_{\sigma \in G}$, where $m: H\to \widetilde{H}$ is a monoidal natural transformation and $m_\sigma:\theta_\sigma\to \widetilde{\theta}_\sigma$ are morphisms in $\C$. The previous data are subject to the following axioms:  $m_e=\id_I$ and the diagrams
$$
\begin{diagram}
\node{\theta_\sigma\widehat{\sigma}(\theta_\tau)U_{\sigma,\tau}'}\arrow{s,l}{\Pi_{\sigma,\tau}} \arrow{e,t}{m_\sigma\otimes \widehat{\sigma}(m_\tau)\otimes\id_{U'_{\sigma,\tau}}} \node{\widetilde{\theta}_\sigma\widehat{\sigma}(\widetilde{\theta}_\tau)U_{\sigma,\tau}'} \arrow{s,r}{\widetilde{\Pi}_{\sigma,\tau}} \\
\node{H(U_{\sigma,\tau})\theta_{\sigma\tau}}\arrow{e,t}{m_{U_{\sigma,\tau}}\otimes m_{\sigma\tau}} \node{\widetilde{H}(U_{\sigma,\tau})\widetilde{\theta}_{\sigma\tau}}
\end{diagram}
\begin{diagram}
\node{\theta_\sigma\widehat{\sigma}(H(V))} \arrow{s,l}{\theta^\sigma_V}\arrow{e,t}{m_\sigma\otimes\widehat{\sigma}(m_V)} \node{\widetilde{\theta}_\sigma\widehat{\sigma}(\widetilde{H}(V))} \arrow{s,r}{\widetilde{\theta}^\sigma_V}\\
\node{H(\widetilde{\sigma}(V))\theta_\sigma} \arrow{e,t}{m_{\widetilde{\sigma}(V)}\otimes m_\sigma} \node{\widetilde{H}(\widetilde{\sigma}(V))\widetilde{\theta}_\sigma}
\end{diagram}
$$commute for all $\sigma,\tau \in G, V\in \C$ (where tensor symbols among objects have been omitted).
\begin{teor}\label{main theorem}
There is a biequivalence between the 2-category of coherent outer $G$-actions and the 2-category of $G$-crossed product tensor categories.
\end{teor}
\begin{proof}
The bijective correspondence between $G$-crossed product tensor categories and coherent outer $G$-action was described in the Section \ref{section action y correspondecia}.

If $T=(H,\theta,\Pi)$ is a 1-cell between coherent outer $G$-action $(\{ \widetilde{\sigma}\}_{\sigma\in G}\chi,\omega)$ and $(\{\widehat{\sigma}\}_{\sigma \in G},\chi',\omega')$ over $\C$ and  $\D$ respectively, then we define an monoidal functor $(T,\psi^T): \C\rtimes G \to \D\rtimes G$ as
\begin{itemize}
  \item $T([V,\sigma])= [H(V)\otimes \theta_\sigma,\sigma], T([f,\sigma])=[H(f)\otimes \id_{\theta_\sigma},\sigma]$ for all $V\in \C, \sigma\in G$,
  \item $\psi^T:T([V,\sigma])\cdot T([W,\sigma])\to T([V,\sigma]\cdot[W,\tau]),$
\end{itemize}where
$$
\begin{diagram}
\node{T([V,\sigma])\cdot T([W,\sigma])} \arrow[4]{s,l}{\psi^T_{[V,\sigma],[W,\tau]}} \arrow{e,t,=}{}
\node{H(V)\theta_\sigma\widehat{\sigma}(H(W)\theta_\tau)U'_{\sigma,\tau}} \arrow{s,r}{\id_{H(V)\theta_\sigma}\otimes\psi^{\widehat{\sigma}}_{H(W),\theta_\tau}\otimes \id_{U'_{\sigma,\tau}}}\\
\node{} \node{H(V)\theta_\sigma\widehat{\sigma}(H(W))\widehat{\sigma}(\theta_\tau)U'_{\sigma,\tau}} \arrow{s,r}{\id_{H(V)}\otimes\theta^{\sigma}_W\otimes\id_{\widehat{\sigma}(\theta_\tau)U'_{\sigma,\tau}}}\\
\node{} \node{H(V)H(\widetilde{\sigma}(W))\theta_\sigma\widehat{\sigma}(\theta_\tau)U'_{\sigma,\tau}}\arrow{s,r}{\psi^H_{V,\widetilde{\sigma}(W)}\otimes \Pi_{\sigma,\tau}}\\
\node{}\node{H(V\widetilde{\sigma}(W))H(U_{\sigma,\tau})\theta_{\sigma\tau}} \arrow{s,r}{\psi^H_{V\widetilde{\sigma}(W),U_{\sigma,\tau}}\otimes \id_{\theta_{\sigma\tau}}}\\
\node{T([V,\sigma]\cdot[W,\tau])}\arrow{e,t,=}{} \node{H(V\widetilde{\sigma}(W)U_{\sigma,\tau})\theta_{\sigma\tau}}
\end{diagram}
$$ (where tensor symbols among objects of $\C$ have been omitted as a space-saving measure). Conversely, given a graded monoidal functor $(T,\psi^T):\C\rtimes G\to \D\rtimes G$, we define a 1-cell $(H,\theta,\Pi)$ as $[H(V),e]=T([V,e])$, $[\theta_\sigma,\sigma]=T([I,\sigma])$, $[\theta^\sigma_V,\sigma]=\psi^T_{[I,\sigma],[V,e]}$, $[\Pi_{\sigma,\tau},\sigma\tau]=\psi^T_{[I,\sigma],[I,\tau]}$.
\smallbreak
Given a 2-cell $\{m_\sigma, m\}_{\sigma \in G}$ between 1-cells $T=(H,\theta,\Pi)$ and $T'=(H',\theta',\Pi')$, we define a monoidal natural isomorphism $m:T\to T'$ between the associated monoidal functors by $m_{[V,\sigma]}=[m_V\otimes m_\sigma,\sigma]$. Conversely, given a monoidal natural isomorphism  $m:T\to T'$, we define a 2-cell by  $[m_V,e]= m_{[V,e]}$, $[m_\sigma,\sigma]=m_{[I,\sigma]}$.
\smallbreak
Finally, in order to see that the 2-categories are biequivalent, note that every crossed product tensor category is equivalent to one of the form $\C_e\rtimes G$. So, every functor between  $\C\rtimes G$ and  $\D\rtimes G$ is monoidally equivalent to one induced by a 1-cell  of the coherent outer $G$-action 2-category, and every monoidal natural transformation is equal to one induced by a 2-cell.
\end{proof}

\section{Braided crossed product tensor categories}\label{seccion 6}

Recall that a braiding for a monoidal category $(\C,\otimes ,I, \alpha)$ is a natural isomorphism $c: \otimes\to \otimes \tau$, where $\tau: \C\times \C\to \C\times \C$ is the flip, and the hexagons

\begin{equation}\label{H1}
\begin{diagram}
\node{}\node{(U\otimes V)\otimes W}\arrow{e,t}{c_{U\otimes V, W}}\node{W\otimes (U\otimes V)}\arrow{se,l}{a_{W,U,V}} \\
\node{U\otimes (V\otimes W)} \arrow{se,r}{\id_{U}\otimes c_{V,W}}\arrow{ne,l}{a_{U,V,W}}\node{}\node{}\node{(W\otimes U)\otimes V}\\
\node{}\node{U\otimes (W\otimes V)}\arrow{e,t}{a_{U,W,V}}\node{(U\otimes W)\otimes V}\arrow{ne,r}{c_{U,W}\otimes \id_V}
\end{diagram}
\end{equation}

\begin{equation}\label{H2}
\begin{diagram}
\node{}\node{U\otimes (V\otimes W)}\arrow{e,t}{c_{U,V\otimes W}}\node{(V\otimes W)\otimes U}\arrow{se,l}{a_{V,W,U}^{-1}} \\
\node{(U\otimes V)\otimes W} \arrow{se,r}{ c_{U,V}\otimes \id_W}\arrow{ne,l}{a_{U,V,W}^{-1}}\node{}\node{}\node{V\otimes (W\otimes U)}\\
\node{}\node{(V\otimes U)\otimes W}\arrow{e,t}{a_{V,U,W}^{-1}}\node{V \otimes (U\otimes W)}\arrow{ne,r}{ \id_V \otimes c_{U,W}}
\end{diagram}
\end{equation}
commute for all $U, V,W\in \C$.
\smallbreak
If a $G$-crossed product tensor category admits a braiding, the group $G$ must be abelian. So, from now on we shall suppose that $G$ is abelian.
\smallbreak
Let $\C$ be a tensor category with a coherent outer $G$-action, such that the tensor category $\C\rtimes G$ admits a braiding $c$. The braiding  $c_{[V],[\sigma]}:[V,\sigma]\to [\sigma_*(V),\sigma]$ defines  natural isomorphisms $c_{V,\sigma}:V\to \sigma_*(V)$. The commutativity of the hexagon \eqref{H1} is equivalent to $c_{V,\sigma}$ is a monoidal natural isomorphism from $\id_\C$ to $\sigma_*$. For that reason, if $\C\rtimes G$ has a braiding, we can suppose that $\sigma_*=\id_\C$ for all $\sigma\in G$.

\begin{defin}
A coherent outer $G$-action shall be called central if $\sigma_*=\id_\C$ for all $\sigma \in G$.
\end{defin}
\begin{obs}
For a central coherent outer $G$-action, the pseudonatural transformations $(\chi_{\sigma,\tau},U_{\sigma,\tau})$ are just elements in $\mathcal{Z}(\C)$ (the center of $\C$), and the modifications $\omega$ are morphisms in $\mathcal{Z}(\C)$.
\end{obs}

\begin{defin}
Let $(\C,c)$ be a braided tensor category, and let $G$ be an abelian group. A braiding for a central  coherent $G$-action over $\C$ is a triple $(\theta^\sigma, \overline{\theta}^\sigma, t_{\sigma,\tau})_{\sigma,\tau\in G}$, where

\begin{itemize}
  \item $\theta^\sigma, \overline{\theta}^\sigma:\id_\C\to \id_\C$ are  monoidal natural isomorphisms,
  \item $t_{\sigma,\tau}: \chi_{\sigma,\tau}\to \chi_{\tau,\sigma}$ are isomorphisms in $\mathcal{Z}(\C)$ for all $\sigma,\tau \in G$,
\end{itemize}
such that $\theta^e =\overline{\theta}^e= \id$, $\theta_I^\sigma =\id_I$, $t_{\sigma,e}=t_{e,\sigma}=\id_I$, and the diagrams

\begin{equation}\label{B1}
\begin{diagram}
\node{Z\otimes U_{\sigma,\tau}}\arrow{se,r}{((\overline{\theta}_Z^{\sigma\tau})^{-1}\overline{\theta}_Z^{\sigma}\overline{\theta}_Z^{\tau})\otimes \id_{U_{\sigma,\tau}}} \arrow{e,t}{\chi_{\sigma,\tau}(Z)}\node{U_{\sigma,\tau}\otimes Z} \arrow{s,r}{c_{U_{\sigma,\tau},Z}} \\
\node{}\node{Z\otimes U_{\sigma,\tau}}
\end{diagram}
\end{equation}

\begin{equation}\label{B2}
\begin{diagram}
\node{Z\otimes U_{\sigma,\tau}}\arrow{se,r}{((\theta_Z^{\sigma\tau})^{-1}\theta_Z^{\sigma}\theta_Z^{\tau})\otimes \id_{U_{\sigma,\tau}}} \arrow{e,t}{c_{U_{\sigma,\tau},Z}}\node{U_{\sigma,\tau}\otimes Z} \arrow{s,r}{\chi_{\sigma,\tau}(Z)} \\
\node{}\node{Z\otimes U_{\sigma,\tau}}
\end{diagram}
\end{equation}

\begin{equation}\label{B3}
\begin{diagram}
\node{}\node{U_{\sigma,\tau}U_{\sigma\tau,\rho}}\arrow[2]{e,t}{\theta_{U_{\sigma,\tau}}^\rho\otimes t_{\sigma\tau,\rho}} \node{}\node{U_{\sigma,\tau}U_{\rho,\sigma\tau}} \arrow{se,l}{w_{\rho,\sigma,\tau}}\\
\node{U_{\tau,\rho}U_{\sigma,\tau\rho}}\arrow{ne,l}{\omega_{\sigma,\tau,\rho}} \arrow{se,r}{t_{\tau,\rho}\otimes \id_{U_{\sigma,\rho\tau}}} \node{}\node{}\node{}\node{U_{\rho,\sigma}U_{\sigma\rho,\tau}}\\
\node{} \node{U_{\rho,\tau}U_{\sigma,\rho\tau}} \arrow[2]{e,t}{\omega_{\sigma,\rho,\tau}}\node{} \node{U_{\sigma,\rho}U_{\sigma\rho,\tau}} \arrow{ne,r}{t_{\sigma,\tau}\otimes \id_{U_{\sigma\rho,\tau}}}
\end{diagram}
\end{equation}

\begin{equation}\label{B4}
\begin{diagram}
\node{}\node{U_{\tau,\rho}U_{\sigma,\tau\rho}}\arrow[2]{e,t}{\overline{\theta}^\sigma_{U_{\tau,\rho}}\otimes t_{\sigma,\tau\rho}}
\node{} \node{U_{\tau,\rho}U_{\tau\rho,\sigma}} \arrow{se,l}{\omega_{\tau,\rho,\sigma}^{-1}}\\
\node{U_{\sigma,\tau}U_{\sigma\tau,\rho}}\arrow{ne,l}{\omega^{-1}_{\sigma,\tau,\rho}} \arrow{se,r}{t_{\sigma,\tau}\otimes\id_{U_{\sigma\tau,\rho}}} \node{}\node{}\node{}\node{U_{\rho,\sigma}U_{\tau,\rho\sigma}}\\
\node{}\node{U_{\tau,\sigma}U_{\tau\sigma,\rho}}\arrow[2]{e,t}{\omega_{\tau,\sigma,\rho}^{-1}}\node{}\node{U_{\sigma,\rho}U_{\tau,\sigma\rho}} \arrow{ne,r}{t_{\sigma,\rho}\otimes \id_{U_{\tau,\sigma\rho}}}
\end{diagram}
\end{equation} commute for all $\sigma,\tau, \rho \in G, Z\in \C$ (where tensor symbols among objects have been omitted as a space-saving measure).
\end{defin}

\begin{teor}\label{teor braid}
Let $(\C,c)$ be a braided tensor category with a coherent central outer $G$-action. Then, there is a bijective correspondence between braidings over $\C\rtimes G$ and braidings over the central coherent outer $G$-action of $\C$.
\end{teor}

\begin{proof}
Let $(\theta,\overline{\theta},t)$ be a braiding for a central coherent outer $G$-action $(\chi,\omega)$. Then, we define a braiding over $\C\rtimes G$ by
\begin{align*}
c_{[V,\sigma],[W,\tau]} &= (c_{V,W}\circ (\theta^\tau_V\otimes\overline{\theta}^\sigma_W))\otimes t_{\sigma,\tau}\\
&= ((\overline{\theta}^\sigma_W\otimes \theta^\tau_V)\circ c_{V,W} )\otimes t_{\sigma,\tau}.
\end{align*}
Conversely, given a braiding $c$ over $\C\rtimes G$, we define a braiding for the coherent central outer $G$-action by
\begin{align*}
    [\theta_V^\sigma,\sigma] := c_{[V],[\sigma]}, \  \  \  [\overline{\theta}_V^\sigma,\sigma]:= c_{[\sigma],[V]},\  \  \text{and}\ \  \ [t_{\sigma,\tau},\sigma\tau]:= c_{[\sigma],[\tau]}.
\end{align*}

Let we denote by $H(U,V,W)$ and $H'(U,V,W)$ the hexagons \eqref{H1} and \eqref{H2}, respectively. Let $\theta^\sigma,\overline{\theta}^\sigma :\id_\C\to \id_\C$ be natural isomorphisms for each $\sigma\in G$, and let $t_{\sigma,\tau}:U_{\sigma,\tau}\to U_{\tau,\sigma}$ be isomorphisms in $\C$. If we set the  following definitions of natural isomorphisms
\begin{align*}
      c_{[V],[\sigma]}:=[\theta_V^\sigma,\sigma], \  \  \  c_{[\sigma],[V]}:=[\overline{\theta}_V^\sigma,\sigma],\  \  \ c_{[\sigma],[\tau]}:=[t_{\sigma,\tau},\sigma\tau],
\end{align*}
it is easy to see that the commutativity of $H([V],[W],\sigma)$ and  $H'([\sigma], [V],[W])$ is equivalent to $\theta^\sigma$ and $\overline{\theta}^\sigma$ be monoidal natural isomorphisms, respectively. The commutativity of $H'([\sigma],[\tau],[Z])$ and $H([\sigma],[V],[\tau])$  is equivalent to $t_{\sigma,\tau}$ be a morphism in $\mathcal{Z}(\C)$. The commutativity of $H([\sigma],[\tau],[Z])$ and $H([Z],[\sigma],[\tau])$  is equivalent to the commutativity of \eqref{B1} and \eqref{B2}, respectively. The commutativity of $H([\sigma],[\tau],[\rho])$ and $H'([\sigma],[\tau],[\rho])$ is equivalent to the  commutativity  of \eqref{B3}, \eqref{B4}, respectively.

\end{proof}

\section{Crossed product tensor categories as quasi-trivial extensions}\label{seccion homotopy}

In \cite{ENO3} Etingof, Nikshych, and Ostrik study fusion categories graded by a finite group, using invertible bimodule categories over fusion categories. They reduce the classification problem of fusion categories graded by a group $G$ to classification (up to homotopy) of maps from $BG$ to classifying spaces of certain higher groupoids. In  this section we shall explain briefly the connection of our results with some results in \cite{ENO3}.

We freely use the language and basic theory of module categories and
tensor product over them, \cite{ENO}, \cite{ENO3}.

In \cite{ENO3} they show that a graded fusion category $\C=\bigoplus_{\sigma\in G}$ determines and it is determined by the following data:

\begin{itemize}
  \item a fusion category $\C_e$, a collection of invertible $\C_e$-bimodule categories $\C_\sigma, \sigma\in G$,
  \item a collection of $\C_e$-bimodule isomorphisms $M_{\sigma,\tau}:\C_\sigma\boxtimes_{\C_e} \C_\tau\to \C_{\sigma\tau}$,
  \item natural isomorphisms of $\C_e$-bimodule functors \[\alpha_{\sigma,\tau,\rho}:M_{\sigma,\tau \rho}(\text{Id}_{\C_\sigma}\boxtimes_{\C_e} M_{\tau,\rho})\to M_{\sigma \tau,\rho}(M_{\sigma,\tau}\boxtimes_{\C_e} \text{Id}_{\C_\rho})\] satisfying the identity
\end{itemize}
\begin{multline}\label{equacion datos graduada}
    M_{\sigma,\tau \rho k}(\id_\sigma\boxtimes_{\C_e}\alpha_{\tau,\rho,k})\circ \alpha_{\sigma,\tau \rho,k}(\text{Id}_{\C_\sigma}
\boxtimes_{\C_e}M_{\tau,\rho}\boxtimes_{\C_e}\text{Id}_{\C_k} )\\
= \alpha_{\sigma,\tau,\rho k}(\text{Id}_{\C_\sigma}
\boxtimes_{\C_e} \text{Id}_{\C_\tau}
\boxtimes_{\C_e} M_{\rho,k}) \circ \alpha_{\sigma\tau,\rho,k}(M_{\sigma,\tau} \boxtimes_{\C_e} \text{Id}_{\C_\rho}
\boxtimes_{\C_e} \text{Id}_{\C_k}),
\end{multline}
for all $\sigma, \tau, \rho, k  \in  G$, where we use the notation Id for the identity functor, and id for the identity morphism.

Following \cite{ENO3} we shall  say that an invertible $\C$-bimodule category $\M$ is quasi-trivial if it is equivalent to $\C$ as a left module category. It is easy to see that if $\M$ is quasi-trivial, then there exists a tensor autoequivalence $\sigma : \C \to \C$, such that $\M= \C$ with the left action of $\C$ by left multiplication, and the right action of $\C$ by right multiplication twisted by $\sigma$.

Given a tensor functor $\sigma:\C\to \C$ we shall denote by $\C^\sigma$ the quasi-trivial $\C$-bimodule associated.

\begin{lem}\label{lem composition}
Let $\C$ be a tensor category and $\sigma, \tau: \C\to \C$ tensor functors. Then  there is a bijective correspondence between $\C$-bimodule functors from $\C^\sigma$ to $\C^\tau$ and pseudonatural transformation from $\sigma$ to $\tau$, and $\C$-bimodule natural transformation and modifications between the pseudonatural transformations associated. Moreover, the tensor product $\C^\sigma\boxtimes_{\C}\C^\tau$ exists and it is equivalent to $\C^{ \sigma\circ \tau}$.
\end{lem}
\begin{proof}
Let $(X,\theta): \sigma \to \tau$ be a pseudonatural transformation. The endofunctor $F_X(-)=(-)\otimes X:\C\to \C$ with the natural isomorphisms $\id_V\otimes\theta_W: F_X(V\otimes^\sigma W)\to F_X(V)\otimes^\tau W$ for all $V,W\in \C$ is $\C$-bimodule functor.

Conversely, suppose that $F:\C^\sigma\to \C^\tau$ is a $\C$-bimodule functor. Let $X=F(1)$, using the natural isomorphisms $F(V)=F(V\otimes 1)\cong V\otimes X$, we can suppose that $F(V)=V\otimes X$ for all $X\in \C$. The natural isomorphisms \[\sigma(W)\otimes X= F(1\otimes^\sigma W) \cong F(1)\otimes^\tau W = X\otimes \tau(W)\] define an pseudonatural transformation.

Let $(X,\theta), (X',\theta'):\C\to \C$ be pseudonatural transformations, and $\alpha: F_X\to F_{X'}$ a $\C$-bimodule natural transformation, where $F_X$ and $F_{X'}$ are the $\C$-bimodule functors associated to $(X,\theta), (X',\theta')$. The morphism $\alpha_{1}:F_X(1)=X\to F_{X'}(1)=X'$ is a modification, and conversely if $\omega: X\to X'$ is a modification, then the natural transformations $\alpha_{V}=\id_V\otimes\omega$, $V\in \C$ is a $\C$-bimodule natural transformation.

For the second part, the category $\C$ with the $\C$-balanced functor $$B_{\C^\sigma,\C^\tau}:\C_\sigma\times \C_\tau \to \C, V\times W\mapsto V\otimes\sigma(W)$$ is a tensor product over $\C$, and it is easy to see that it is equal to $\C^{\sigma\circ \tau}$.
\end{proof}

\begin{prop}
A coherent outer $G$-action $(\{ \sigma_*\}_{\sigma\in G},\chi,\omega)$ defines  data $(\C^\sigma,$ $ M_{\sigma,\tau},\alpha )_{\sigma,\tau\in G}$  that satisfy the equation \eqref{equacion datos graduada} and conversely a data $(\C^\sigma, M_{\sigma,\tau},\alpha )_{\sigma,\tau\in G}$ that satisfy the equation \eqref{equacion datos graduada} and where $\C^\sigma =\C$ as left $\C$-module categories for all $\sigma \in G$, defines a coherent outer action.
\end{prop}
\begin{proof}
Using the Lemma \ref{lem composition} it is easy to see that the composition and the tensor product of pseudonatural transformations and modifications in Bieq$(\C)$ correspond to the composition and the tensor product of the $\C$-bimodule functor and natural transformations associated.

Now, if $(\{ \sigma_*\}_{\sigma\in G},\chi,\omega)$ is a coherent outer $G$-action $\chi$ defines an equivalence of  $\C$-bimodule $M_{\sigma,\tau}:\C_\sigma\boxtimes_{\C_e} \C_\tau\to \C_{\sigma\tau}$, and the modifications $\omega_{\sigma,\tau,\rho}$ define $\C$-bimodule natural isomorphism $\alpha_{\sigma,\tau,\rho}:M_{\sigma,\tau \rho}(\text{Id}_{\C_\sigma}\boxtimes_{\C_e} M_{\tau,\rho})\to M_{\sigma \tau,\rho}(M_{\sigma,\tau}\boxtimes_{\C_e} \text{Id}_{\C_\rho})$, it is a straightforward verification that the equation \eqref{equacion datos graduada} is equivalent to commutativity of the diagram \eqref{coherencia}.
\end{proof}

\bibliographystyle{amsalpha}

\end{document}